\documentclass[11pt]{article}
 
\author{
Barbara  Niethammer
\\
Inst.\ f\"ur Angew.\ Math.\\
Universit\"at Bonn\\
Wegelerstr.\ 6\\
53115 Bonn, Germany\\
\and
Robert L. Pego\\
Dept.\ of Mathematics \& \\
Inst.\ Phys.\ Sci.\ Tech.\\
University of Maryland\\
College Park, MD 20742 USA
}
\date{April 1998}

\usepackage{amsfonts}
\newcommand{\nwc}{\newcommand}

\newtheorem{prop}{Proposition}[section]
\newtheorem{lemma}[prop]{Lemma}

\newtheorem{theorem}[prop]{Theorem}
\newtheorem{corollary}[prop]{Corollary}

\newenvironment
  {proof}{\noindent {\bf Proof:}}{\bigskip}
 
\newcommand {\eps}{\varepsilon}
\nwc{\D}{\partial}
\nwc{\vdot}{\Lambda}
\nwc{\F}{\varphi}

\newcommand {\R} {{\mathbb{R}}}
\nwc{\third}{{1/3}}
\nwc{\RCD}{rcd([0,1])}
\nwc{\BDD}{bdd([0,1])}
\nwc{\Psp}{{\cal P}_0}  
\nwc{\loc}{{\rm \scriptstyle loc}}
\nwc{\gi}{^{\dag}}
\nwc{\git}{^{\dag\dag}}
\nwc{\vnu}{{\hat v}}
\nwc{\nuv}{{\hat\nu}}
\nwc{\nuone}{\nu^{(1)}}
\nwc{\nutwo}{\nu^{(2)}}

\nwc{\vtup}{\overline{\D}_tv}
\nwc{\vtdn}{\underline{\D}_tv}
\nwc{\delvh}{\delta^hv}
\nwc{\delvhj}{\delta^{h_j}v}
\nwc{\Ccr}{{C_c(\R^+)}}
\nwc{\dist}{\mathop{\rm dist}\nolimits}
\nwc{\ess}{\mathop{{\rm ess}}\nolimits}
\nwc{\esssup}{\mathop{\rm ess\,sup}}

\renewcommand{\AE}{{a.e.\ }} 
\nwc{\ds}{\displaystyle}
\nwc{\eqref}[1]{(\ref{#1})}
\newcommand {\bedis} {\begin{displaymath}}
\newcommand {\edis} {\end{displaymath}}
\newcommand{\newbeqna} {\renewcommand {\arraystretch} {2}
                        \begin {displaymath} \begin {array}{crcl}}
\newcommand{\neweqna}{\end{array} \end {displaymath}}
 
\newcommand{\fbeqna}{\renewcommand {\arraystretch} {1.3}
\begin {displaymath}\begin{array}{rcll}}
\newcommand{\feqna}{\end{array}\end{displaymath}}
 
\newcommand {\beqna} {\begin{eqnarray*}}
\newcommand {\eqna} {\end{eqnarray*}}
\newcommand {\beqn} {\begin{eqnarray}}
\newcommand {\eqn} {\end{eqnarray}}
\newcommand {\be} {\begin{equation}}
\newcommand {\ee}{\end{equation}}

\begin{document}
\title{ On the initial-value problem in the \\
Lifshitz-Slyozov-Wagner theory of Ostwald ripening}
\maketitle
 
\begin{abstract}
The LSW theory of Ostwald ripening concerns the time evolution 
of the size distribution of a dilute system of particles
that evolve by diffusional mass transfer with a common mean field.
We prove global existence, uniqueness and continuous dependence
on initial data for measure-valued solutions with compact support 
in particle size. 
  These results are established with respect to a natural topology on the
space of size distributions, one given by the Wasserstein
metric which measures the smallest maximum volume change required to
rearrange one distribution into another.

\end{abstract}

\section{Introduction}

The classical theory of Ostwald ripening, formulated by Lifshitz and 
Slyozov \cite{LS} and Wagner \cite{Wa} concerns the evolution of the 
size distribution of a large number of small particles of one phase
embedded in a matrix of another phase. Particles are assumed to be
widely separated spheres that evolve by diffusional mass transfer with
a common mean field. In the late stages of the phase transformation,
diffusion is quasi-steady and the particle growth rate is determined
by the mass flux at the particle boundary. The mass flux is proportional
to the gradient of a potential that is harmonic, is proportional to curvature on the 
particle boundaries, and is close to constant in the mean field
between particles.

In appropriate units, it is found that any particle
radius $R(t)$ evolves according to 
\be
\frac{dR}{dt} = V(R,R_c(t))\colon= 
\frac{a}{R^2}\left(\frac{R}{R_c(t)}-1\right),
\label{1Rdot}
\ee
where $a$ is a constant and the critical radius $R_c(t)$ is the same 
for all particles. The value of $R_c(t)$ is determined from conservation
of mass. If mass changes in the diffusion field can be neglected, 
the particle
volume is conserved and one finds that the critical radius equals 
the average radius of currently existing particles. Particles with
radius larger than $R_c(t)$ are growing, and particles with smaller
radius shrink and can disappear in finite time. 

Classically, the size distribution of particles is described by a
particle radius distribution $n(t,R)$. This is a normalized
number density that we may scale so that $\int_0^R n(t,r)\,dr$ is 
the number of (currently existing) particles with radius 
less than $R$, divided by the number $N$ of 
initially existing particles. The number of particles with
size between $R_1(t)$ and $R_2(t)$ for any two solutions of
\eqref{1Rdot} is conserved, so $n(t,R)$ should satisfy the conservation law
\be
\D_t n +\D_R(Vn)=0,
\label{1neqn}
\ee
where the critical radius is given by 
\be
R_c(t) = \int_0^\infty Rn(t,R)\,dR\left/\int_0^\infty n(t,R)\,dR\right. .
\label{1Rc}
\ee
The initial number density $n_0(R)=n(0,R)$ 
satisfies $\int_0^\infty n_0(R)\,dR=1$ in this normalization. 

Our aim in this paper is to develop a satisfactory
theory of well-posedness for the initial value problem for the 
particle size distribution. From the physical point of view, 
it is reasonable to suppose that a positive fraction of the particles 
can have the same radius, in which case the size distribution
contains one or more Dirac deltas. Mathematically, the ideal
is to allow the initial data
$n_0(R)\,dR$ to be an arbitrary probability measure such that
the total volume $
\int_0^\infty \frac{4}3\pi R^3n_0(R)\,dR
$
is finite. 

It will be convenient to
work with particle volume $v$ instead of radius $R$, and to
work with a cumulative number distribution function $\F$ instead
of the number density $n$. We say that
\be
\mbox{$\F$ is the fraction of (initially existing)
particles with volume $\ge v$.}
\label{1Fdef}
\ee
As a function of volume $v$ at time $t$, $\F(t,v)$ is a monotonically
decreasing function which is left continuous at jumps with $\F(t,0)=1,$
and $\int_0^\infty \F(t,v)\,dv$ (the total volume) is independent of time.
The particle volume distribution, defined by $f(t,v)\,dv=-d\F(t,v)$
for each fixed $t$, is formally related to $n$ via $f(t,v)\,dv=n(t,R)\,dR$. 

We normalize the time scale by the factor $4\pi a$ and let 
$\theta(t)=(4\pi R_c(t)^3/3)^{-1/3}$,
so that the volume $v(t)$ of any existing particle should satisfy
\be
\frac{dv}{dt} = \vdot(v,\theta(t)) := v^\third\theta(t)-1.
\label{1vdot}
\ee
If $v(t)$ is a positive solution of \eqref{1vdot} on some time 
interval, then $\F(t,v(t))$ should remain constant. This means
$\F(t,v)$ should be a solution of the hyperbolic equation
\be
\D_t\F+\vdot(v,\theta(t))\D_v \F = 0,
\label{1Feqn}
\ee
whose characteristics satisfy \eqref{1vdot}.
The value of $\theta(t)$ is obtained from $\F$
in terms of Riemann-Stieltjes integrals by 
\be
\theta(t)= \int_{0^+}^\infty d\F(t,v) \left/
{\int_0^\infty v^\third \,d\F(t,v)}\right. .
\label{1theta}
\ee
The numerator is $-1$ times the quantity $\F_0(t):=\lim_{v\to0}\F(t,v)$,
which is the fraction of initially existing particles that still
exist at time $t$.

It turns out to be still  better to regard the volume $v$ as a 
function of the fraction $\F$, $0\le\F\le1$. 
We take the map $\F\mapsto v(t,\F)$ to be right continuous
and decreasing with $v(t,1)=0$. Mathematically, given $\F(t,v)$ we obtain
$v(t,\F)$ via the prescription
\be
v(t,x) = \sup\{y\mid \F(t,y)>x\} \quad \mbox{for $0\le x<1=\max \F$}.
\label{1Finv}
\ee
This is most easily understood when the size distribution
corresponds to a finite number of particles. 
If we list the particle volumes in decreasing order, 
$v_0(t)\ge\ldots\ge v_{N-1}(t)$, then $v(t,\F)=v_j$ for
$\F\in[j/N,(j+1)/N)$. We shall call $\F\mapsto v(t,\F)$ a 
{\it volume ordering} for the system at time $t$.

For technical simplicity we shall assume that the particle volumes in
the system are bounded. This seems reasonable physically, and corresponds
to assuming that the particle volume distribution has compact support 
in $v$. We then introduce function spaces as follows.
Let $\RCD$ be the set of functions $v\colon[0,1]\to\R$
that are right continuous, decreasing, and satisfy $v(1)=0$. 
(To be precise, we say $v$ is {\it decreasing} 
if $v(x_1)\le v(x_2)$ whenever $x_1\ge x_2$, and similarly for 
{\it increasing}. A decreasing function need not be strictly decreasing.)
The set $\RCD$ is contained in the space $\BDD$ of real-valued
bounded functions on $[0,1]$, equipped with the sup norm 
$\|v\| = \sup_\F |v(\F)|$.  $\RCD$ is a 
complete metric space in the induced topology.

If $X$ is a Banach space and $I\subset\R$ is an interval, then
$C(I,X)$ is the space of continuous $X$-valued functions on $I$,
and $L^\infty_\loc(I)$ is the space of equivalence
classes of measurable functions locally bounded on $I$,
where two functions are considered equivalent if they agree almost
everywhere. 

Our main results are the following. 

\begin{theorem} \label{T.eu}
(Global existence and uniqueness)
Let $v_0\in\RCD$. Then there exist unique functions 
$\theta\in L^\infty_\loc(0,\infty)$ and 
$v\in C([0,\infty),\RCD)$, such that
\[
\int_0^1 v(t,\F)\,d\F = \int_0^1 v_0(\F)\,d\F
\]
for all $t\ge0$, and
\[
v(t,\F)= v_0(\F)+\int_0^t (v(s,\F)^\third\theta(s)-1)\,ds
\]
for all $(t,\F)$ such that $v(t,\F)>0$.
\end{theorem}

\begin{theorem} \label{T.lip} 
(Continuous dependence on initial data)
Given positive constants $T$ and $C_0$, there exists a positive constant 
$C$ such that,
if $(v_1,\theta_1)$ and $(v_2,\theta_2)$ are two solutions with
the properties stated for $(v,\theta)$ in Theorem \ref{T.eu},
and if $\max(v_1(0,0),v_2(0,0))\le C_0$, 
then
\[
\sup_{
{0\le t\le T}
}
\|v_1(t,\cdot)-v_2(t,\cdot)\| \le C \|v_1(0,\cdot)-v_2(0,\cdot)\|.
\]
Consequently, the map $v_0\mapsto v$ is locally Lipschitz 
from $\RCD$ into $C([0,T],\RCD)$.
\end{theorem}

Our strategy to prove these results 
at the same time justifies a method of numerical approximation 
for the problem that has a direct physical interpretation. 
We first consider solutions that are piecewise
constant, taking a finite number of values $v_0(t)>\ldots>v_{N-1}(t)$,
as is the case for a finite number of particles.
We show that these solutions are determined on a succession of time
intervals by solving finite systems of coupled ordinary differential 
equations with a number of components that decreases as the smallest
particles vanish. Once we prove the continuity estimate in Theorem
\ref{T.lip} (at first for initial data near to each other), 
uniqueness is immediate and existence for general
initial data in $\RCD$ follows by an approximation argument.

\smallskip
The solutions constructed in Theorem \ref{T.eu} correspond to 
measure-valued weak solutions of the evolution equation 
\be
\D_tf+\D_v(\vdot(v,\theta(t))f) = 0
\label{1fevol}
\ee
for the particle volume distribution. This means that at each
time $t$, the formal expression $f(t,v)\,dv$ corresponds to a
probability measure $\nu_t$ having compact support in $[0,\infty)$,
the set of volumes.
The notion of distance used in Theorem \ref{T.lip} has an 
interpretation as a natural metric on the space $\Psp$ of such
probability measures. This metric measures the smallest `maximum
volume change' required to rearrange one volume distribution
into another. Mathematically it is the $L^\infty$ Wasserstein metric \cite{GS,Rach},
which we denote by $d_\infty$.
In section 3 we shall establish 
the relationship between $v(t,\F)$ and $\nu_t$, and deduce the
following result as a corollary of Theorems \ref{T.eu} and \ref{T.lip}.

\begin{theorem} \label{T.meas}
Let $\Psp$ denote the set of probability measures on $[0,\infty)$ 
of compact support, with topology given by $d_\infty$, the $L^\infty$
Wasserstein metric. Given $\nu_0\in\Psp$, there exists
a unique
$\theta\in L^\infty_\loc(0,\infty)$ and a unique map $t\mapsto\nu_t$
that is locally Lipschitz from $[0,\infty)$ into $\Psp$, such that
$(\theta,\nu)$ is a volume-conserving weak solution of \eqref{1fevol}, 
in the sense that
\[
\int_0^\infty v\,d\nu_t(v) = \int_0^\infty v\,d\nu_0(v)
\]
for all $t\ge0$, and 
\[
\int_0^\infty\int_0^\infty \D_t\zeta(t,v)+\vdot(v,\theta(t))
\D_v\zeta(t,v)\ d\nu_t(v)\,dt = 0
\]
for all smooth $\zeta\colon(0,\infty)\times(0,\infty)\to\R$ 
with compact support. 

Furthermore, given any $T>0$, $C_0>0$, there exists $C>0$  such that,
if two such weak solutions $(\theta_1,\nuone)$, $(\theta_2,\nutwo)$ are
given, such that the supports of $\nuone_0$ and $\nutwo_0$ are contained
in $[0,C_0]$, then
\[
\sup_{0\le t\le T} d_\infty(\nuone_t,\nutwo_t) \le C\, d_\infty(\nuone_0,\nutwo_0).
\]
\end{theorem}

It is arguably natural from the physical point of view to
measure distance between volume distributions
by using the Wasserstein distance as is done here. 
A physically reasonable notion of distance should reflect 
in a plausible way the effect of small perturbations of the system
on size distributions. In late-stage Ostwald ripening one imagines
that the nucleation or destruction of large particles is unlikely.
Thus the topology should not make it `easy' to change the 
number of large particles. It is plausible, rather, that small perturbations
to the system would involve small changes to particle volumes.
These notions are captured here by the use of the sup norm distance
between volume orderings, and this is equivalent to using the $L^\infty$ 
Wasserstein metric to compare volume distributions.

In section 4 we briefly treat a related, but simpler, case that 
arises in LSW theory, in which mass variations in the diffusion field 
are not neglected. 
In this case it is not total particle volume that is conserved in time,
but rather a quantity of the form
\[
a\theta(t)+\int_0^1v(t,\F)\,d\F,
\]
where $a>0$ is constant.
The evolution of particle volumes is still given by \eqref{1vdot}, but $\theta$ 
is now determined directly from the conserved quantity.

\section{A priori estimates and well-posedness}

In order to prove the a priori estimate stated in Theorem~\ref{T.lip},
we need a pair of lemmas that yield strengthened variants 
of Gronwall's inequality.

\begin{lemma} \label{L.gron}
Suppose $G\colon[0,T]\to\R$ is increasing with $G(0)=0$, $K\ge0$ is
a constant and $f\colon[0,T]\to\R$ is continuous and satisfies
\[
0\le f(t)\le K+\int_{0^+}^{t} f(s)\,dG(s),\quad 0\le t\le T.
\]
Then $f(t)\le Ke^{G(t)}$ for $0\le t\le T$.
\end{lemma}

\begin{proof}
Let
\[ U(t)= K+\int_{0^+}^{t} f(s)\,dG(s), \]
then $U(0)=K$ and $U$ is increasing. To prove the lemma it suffices
to show that $e^{-G}U\le K$. Let $\{t_j\}_{j=0}^n$ be a partition of
$[0,T]$ and define
\[
\Delta t = \sup_{1\le j\le n}(t_j-t_{j-1}), \quad
\epsilon(\Delta t) = \sup_{|t-s|\le\Delta t}|f(t)-f(s)|.
\]
Put $U_j=U(t_j)$, $G_j=G(t_j)$. Then
\begin{eqnarray*}
e^{-G_{j+1}}U_{j+1} - e^{-G_j}U_j &=&
e^{-G_{j+1}}(U_{j+1} - U_j)- U_j(e^{-G_j} - e^{-G_{j+1}})
\\ &=& e^{-G_{j+1}} \left( \int_{t_j}^{t_{j+1}} f(s)\,dG(s) -
U_j\left(e^{G_{j+1}-G_j} -1\right) \right)
\\&\le& e^{-G_{j+1}} \left(f(t_j)+\epsilon(\Delta t) - U_j \right)
(G_{j+1}-G_j) 
\\&\le& \epsilon(\Delta t)(G_{j+1}-G_j), 
\end{eqnarray*}
where we used that $e^x-1\ge x$ for all $x$ and $f(t_j)\le U_j$.
Summing, we find that $e^{-G_j}U_j\le K+\epsilon(\Delta t)G_j$ for all
$j$. Since the partition is arbitrary, $\epsilon(\Delta t)$ can be made
arbitrarily small and the result follows.
\end{proof}

{  
\nwc{\tilt}{{\tilde{t}}}
\nwc{\tilG}{{\tilde{G}}}
\nwc{\Ti}{\tau}
\nwc{\tilTi}{{\tilde{\tau}}}
\nwc{\tpf}{{t+f(t)}}

\begin{lemma} \label{L.diff}
Suppose $G\colon[0,T]\to\R$ is increasing,
and $f\colon[0,T]\to\R$ is continuous and nonnegative 
and increasing.
Then as long as $0\le t+f(t)\le T$ we have
\[
\int_0^t \left(G(s+f(s))- G(s)\right)ds \le 
\int_0^{f(0)}(G(f(0))-G(s))\,ds +\int_0^t f(s)\,d\tilG(s)
\]
where $\tilG(s)= G(s+f(s))$. 
\end{lemma}

\begin{proof}
Let $Q$ denote the quantity on
the right hand side of the desired inequality. 
Observe that since $G$ is increasing, we have that
\[
Q + \int_t^\tpf (G(s+f(s))-G(s))ds \ge 
Q +\int_t^\tpf \tilG(s)\,ds - G(t+f(t))f(t) .
\]
Since $G(t+f(t))=\tilG(t)$, after integrating by parts in the 
Riemann-Stieltjes integral and cancelling boundary terms we find that the 
last right hand side equals
\begin{eqnarray*}
&&-\int_0^{f(0)} G(s)\,ds -\int_0^t \tilG(s)\,df(s)+\int_t^\tpf\tilG(s)\,ds\\ 
&&\quad=\ -\int_0^{f(0)}G(s)\,ds -\int_0^t\tilG(s)\,d(s+f(s)) 
+\int_0^\tpf\tilG(s)\,ds\\
&&\quad=\ \int_0^\tpf \left(G(s+f(s))-G(s)\right)\,ds .
\end{eqnarray*}
Cancelling the part of the integral from $t$ to $t+f(t)$
finishes the proof.
\end{proof}
} 

Next we establish some basic properties of solutions of the initial value
problem as described in Theorem~\ref{T.lip}. Fixing $T>0$, we shall
consider $t\in[0,T]$. Let $\theta\in L^\infty(0,T)$ be positive and let
$v\in C([0,T],\RCD)$ be such that 
\be
\int_0^1 v(t,\F)\,d\F = \int_0^1 v(0,\F)\,d\F
\label{2vint}
\ee
for all $t$ and 
\be
v(t,\F)=v(0,\F)+\int_0^t(v(s,\F)^{1/3}\theta(s)-1)\,ds
\label{2veqn}
\ee
whenever $v(t,\F)>0$. By scaling, we may assume $\int_0^1 v(t,\F)\,d\F=1$
for all $t\in[0,T]$. 

From \eqref{2veqn} it follows that $t\mapsto v(t,\F)$ is Lipschitz and
satisfies
\be
\frac{\D v}{\D t}= v^{1/3}\theta(t) -1
\label{2vode}
\ee
for almost every $t$ in any interval where $v>0$. 
Since $v^{1/3}\theta-1\le-\frac12$ for $v<\eps_0$ where 
\[
\eps_0^{-1}=8 \esssup_{0\le t\le T}\theta(t),
\]
it follows easily that if $v(t_0,\F)=0$ then $v(t,\F)=0$ for all $t\ge t_0$.

We define $\bar{v}(t)= v(t,0)=\max_\F v(t,\F)$ and
with the notation $a\wedge b=\min(a,b)$ we define
\beqna
\bar t(\F) &=& \inf\{t\in[0,T]\mid v(t,\F) = 0\}\wedge T, \\ 
\bar \F(t) &=& \sup\{\F\in[0,1]\mid v(t,\F)>0\}.
\eqna
The functions $\bar t$ and $\bar \F$ are decreasing functions, and
$\bar \F(t)>0$ for all $t$, since $v(t,\cdot)$ can never vanish
identically by volume conservation. 
We call $\bar t(\F)$ the {\it vanishing time}
for $v(t,\F)$ at $\F$ if $\bar t(\F)<T$ (but note that $\bar t(\F)=T$ if
$v(T,\F)>0$).

\begin{lemma} \label{L.theta1}
For almost every $t\in[0,T]$ we have
\[
0<\theta(t) = \frac{\bar \F(t)}{\ds \int_0^1 v(t,\F)^{1/3}\,d\F} \le
\bar v(t)^{2/3} \le (e^t\bar v(0))^{2/3}.
\]
\end{lemma}

\begin{proof} Evaluate \eqref{2veqn} at $\min(t,\bar t(\F))$ and
integrate over $\F\in[0,1]$. Changing the order of integration and using
the fact that $\bar v(\bar t(\F),\F)=0$ if $\bar t(\F)<t$, we obtain
\begin{eqnarray*}
0&=& \int_0^{\bar \F(t)} v(t,\F)\,d\F - \int_0^1 v(0,\F)\,d\F
\\&=& \int_0^1\int_0^{\min(t,\bar t(\F))} (v(s,\F)^{1/3}\theta(s)-1)\,ds\,d\F
\\&=& \int_0^t\int_0^{\bar \F(s)} (v(s,\F)^{1/3}\theta(s)-1)\,d\F\,ds .
\end{eqnarray*}
Since $t$ is arbitrary the formula for $\theta(t)$ follows.
To get the inequalities, we use that $\bar \F(t)\le1$ and
$\int_0^1 v^{1/3}\,d\F \ge \bar v(t)^{-2/3}\int_0^1 v\,d\F$. Then since
$d\bar v/dt \le \bar v^{1/3}\theta\le \bar v$ we find
\be
\bar v(t) \le e^t\bar v(0).
\label{2vbarest}
\ee
\end{proof}

\begin{lemma} \label{L.smallv}
Whenever $v(t_1,\F)<\eps_0$, 
we have $\D v/\D t<-\frac12$ for
almost every $t\in[t_1,\bar t(\F)]$ and
\[
\frac12 (\bar t(\F)-t)  \le v(t,\F) <\eps_0 - \frac12(t-t_1)
\]
for all $t\in[t_1,\bar t(\F)]$.
\end{lemma}

\begin{proof} $v<\eps_0$ implies $v^{1/3}\theta-1<-\frac12$ almost everywhere, 
and the results follow easily.
\end{proof}

\begin{corollary} \label{C.vest}
 There is a constant $C=C(T,C_0)$ such that
\[
\int_0^{\bar t(\F)} v(t,\F)^{-2/3}\,dt \le C
\]
for all $\F\in[0,1]$. Furthermore, the function $\beta$ given by
\[
\beta(t)=\int_0^{\bar \F(t)} v(t,\F)^{-2/3}\,d\F
\]
is finite for \AE $t\in[0,T]$ and $\int_0^T\beta(t)\,dt \le C$.
\end{corollary}

\begin{proof} The first assertion follows directly from the estimates of
the preceding lemma. The second follows from Fubini's theorem.
\end{proof}

Our plan now is to first prove a restricted version of Theorem~\ref{T.lip}, 
for two solutions that are initially close together.
This restricted result will suffice to establish the existence and uniqueness
theorem, after which the results of Theorem \ref{T.lip} without restriction 
can be proved.

\begin{prop} \label{P.lip}
Given $T>0$, $C_0>0$, there exist $C>0$ and $\delta>0$
such that the bound asserted in Theorem~\ref{T.lip} holds under the additional
assumption that 
\[
\|v_1(0,\cdot)-v_2(0,\cdot)\| \le \delta.
\]
\end{prop}

To start the proof of this restricted version of
Theorem~\ref{T.lip}, we suppose that $T, C_0>0$ 
are given and put
\[
\eps_1=(8e^TC_0)^{-1}.
\]
We suppose that $(\theta_1,v_1)$ and $(\theta_2,v_2)\in L^\infty(0,T)\times
C([0,T],\RCD)$ are two solutions of \eqref{2vint} and \eqref{2veqn} such that
$\max(v_1(0,0),v_2(0,0))\le C_0$.
We define
\[
M(t) = \sup_{0\le s\le t} \|v_1(s,\cdot)-v_2(s,\cdot)\|
\]
and assume that $M(0)<\eps_1$.

\begin{lemma} \label{L.Mest}
There is a constant $C_1=C_1(T,C_0)$ such that for $0\le t\le T$ we have
\[
M(t)\le C_1\left(M(0)+\int_0^t |\theta_1(s)-\theta_2(s)|\,ds \right).
\]
\end{lemma}

\begin{proof}
Fix $\F\in[0,1]$. We suppose that $\bar t_1(\F)\ge \bar t_2(\F)$ without loss
of generality. For $t\in[0,\bar t_2(\F)]$ we may write
\begin{eqnarray*}
v_1(t,\F)-v_2(t,\F) &=& v_1(0,\F)-v_2(0,\F)
+\int_0^t v_2(s,\F)^{1/3}(\theta_1(s)-\theta_2(s))\, ds
\\&& 
+\int_0^t \theta_1(s)(v_1(s,\F)^{1/3}-v_2(s,\F)^{1/3})\,ds .
\end{eqnarray*}
Using the bounds above for $\theta_1$ and $v_2$, and the fact that
$|a-b|\le|a^3-b^3|/a^2$ whenever $a, b>0$, with $C_*=(e^TC_0)^{1/3}$
we obtain the estimate
\begin{eqnarray}
|v_1(t,\F)-v_2(t,\F) | &\le& |v_1(0,\F)-v_2(0,\F)|
 + C_* \int_0^t |\theta_1(s)-\theta_2(s)|\, ds \nonumber
\\&& + C_*^2 \int_0^t v_1(s,\F)^{-2/3}
|v_1(s,\F)-v_2(s,\F) |\,ds. 
\label{2vvest}
\end{eqnarray}

For $t\in[\bar t_2(\F),\bar t_1(\F)]$, we have $v_2(t,\F)=0$ and may write
\[
v_1(t,\F)\le v_1(\bar t_2(\F),\F) + C_*^2\int_{\bar t_2(\F)}^t v_1(s,\F)^\third\,d\F.
\]
Using \eqref{2vvest} with $t=\bar t_2(\F)$ to estimate $v_1(\bar t_2(\F),\F)$,
we find that \eqref{2vvest} is valid for all $t\in[0,\bar t_1(\F)]$.
Gronwall's inequality then yields that 
\beqna
&& 
 \exp\left(-C_*^2 \int_0^t v_1(s,\F)^{-2/3}ds \right) |v_1(t,\F)-v_2(t,\F) |
\\ && \qquad \le\ |v_1(0,\F)-v_2(0,\F)|+ C_* \int_0^t
|\theta_1(s)-\theta_2(s)|\, ds .
\eqna
Using Corollary \ref{C.vest} completes the proof.
\end{proof}

\begin{lemma} \label{L.Fest}
Suppose $M(t)\le\eps_1$ for $0\le t\le\tau$. Then
\[
|\bar \F_1(t)-\bar \F_2(t)| \le
\bar \F_1(t)-\bar \F_1(t+2M(t)) + \bar \F_2(t)-\bar \F_2(t+2M(t)) .
\]
as long as $t+2M(t)\le\tau$.
\end{lemma}

\begin{proof}
Fixing $t$, by relabeling we can assume $\bar \F_1(t)\le \bar \F_2(t)$.
For $\F\in[\bar \F_1(t),\bar \F_2(t)]$, $s\in[t,\tau]$ we have
$v_1(s,\F)=0$ and $v_2(s,\F)\le M(s)\le\eps_1$ by assumption. By
Lemma~\ref{L.smallv}, for $s\le\bar t_2(\F)$ we have $\D v_2/\D t\le-\frac12$
and therefore $\bar t_2(\F)\le \min(t+2M(t),T)$. Hence
$\bar \F_2(t+2M(t))\le \bar \F_1(t)$, and the result follows.
\end{proof}

\begin{lemma} \label{L.thetaest}
There is a constant $C_2=C_2(T,C_0)$ and an increasing function
$H\colon[0,T]\to\R$ depending on $v_1$ and $v_2$, satisfying $H(0)=0$
and $H(T)\le C_2$,
such that if $M(t)\le\eps_1$ for $0\le t\le\tau$, then
\[
\int_0^t |\theta_1(s)-\theta_2(s)|\,ds \le C_2M(0)+ \int_{0^+}^t M(s)\,dH(s)
\]
as long as $t+2M(t)\le\tau$.
\end{lemma}

\begin{proof}
Using that $\int v_j^{1/3}d\F\ge \bar v_j^{-2/3}\ge C_*^{-2}$,
from the formula for $\theta(t)$ we obtain that
\[
|\theta_1(t)-\theta_2(t)| \le 
C_*^2|\bar \F_1(t)-\bar \F_2(t)|
+ C_*^4 \int_0^1|v_1^{1/3}-v_2^{1/3}|\,d\F.
\]
Let $\F_+(t)=\max(\bar \F_1(t),\bar \F_2(t))$, then for $\F<\F_+$ we have
\[
|v_1^{1/3}-v_2^{1/3}| \le \frac{|v_1-v_2|}{v_1^{2/3}+v_2^{2/3}} .
\]
Note that from Corollary \ref{C.vest}, it follows that with 
$t_+(\F)=\max(\bar t_1(\F),\bar t_2(\F))$ we have
\[
\int_0^{t_+(\F)}\frac{1}{v_1^{2/3}+v_2^{2/3}} dt \le C(T,C_0).
\]
By Fubini's theorem it follows that the function defined by 
\[
h_0(t)= \int_0^{\F_+(t)} \frac{1}{v_1^{2/3}+v_2^{2/3}} d\F
\]
is finite for \AE $t$ and is integrable with
$\int_0^T h_0(t)\,dt\le C(T,C_0)$.  Then we have
\be
 \int_0^1 |v_1^{1/3}-v_2^{1/3}|\,d\F \le M(t) h_0(t)
\label{2.v1/3}
\ee
for \AE $t\in[0,T]$.

Next, for $j=1$ and $2$ we invoke Lemma \ref{L.diff} with $G(t)=-\bar\F_j(t)$,
$f(t)=2M(t)$, and conclude that as long as $t+2M(t)\le\tau$, then
\[
\int_0^t \bar\F_j(s)-\bar\F_j(s+2M(s)) \,ds \le 2M(0)+\int_{0^+}^t 2M(s)\, dH_j(s)
\]
where $H_j(t)=-\bar\F_j(t+2M(t))+\bar\F_j(2M(0))$.
Evidently $H_j$ satisfies $H_j(t)\le1$ for all $t$.

Putting these estimates together with the result of Lemma~\ref{L.Fest},
we find that
\[
\int_0^t| \theta_1(s)-\theta_2(s)|\,ds \le 4C_*^2 M(0)+
\int_{0^+}^t M(s)\,dH(s)
\]
where
\[
H(t)=2C_*^2(H_1(t)+H_2(t))+ C_*^4 \int_0^t h_0(s)\,ds.
\]
The desired result follows.
\end{proof}

The proof of Proposition~\ref{P.lip} uses a continuation argument based on
the estimates above together with the estimate
\be
M(\tau)-M(t)\le 2C_*^3(\tau-t)
\label{2Mdiff}
\ee
whenever $0\le t\le\tau\le T$, which follows from $|\D v/\D t|\le C_*^3$.
Since $M$ is increasing, we can find $\tilde T\le T$ such that
$\tilde T+2M(\tilde T)= T$.
With $\tau=t+2M(t)$, inequality \eqref{2Mdiff} yields
\be
M(t+2M(t))\le M(t)(1+4C_*^3)
\label{2Mdelay}
\ee
whenever $t\le\tilde T$.
Now let
\[
\Omega=\{t\in[0,\tilde T]\mid M(t+2M(t))\le\eps_1\}.
\]
If $M(0)\le\delta_0:=\eps_1/(1+4C_*^3)$, then $0\in\Omega$ so
$\Omega$ is nonempty, and clearly $\Omega$ is closed. We claim $\Omega$
is open in $[0,\tilde T]$ if $M(0)$ is sufficiently small.

Given any $t_1\in\Omega$ we can apply Lemmas~\ref{L.Mest} and
\ref{L.thetaest} to deduce that
\be
M(t)\le C_1(1+C_2)M(0)+C_1\int_{0^+}^t M(s)\,dH(s)
\label{2gronM}
\ee
for $0\le t\le t_1$. Then Lemma \ref{L.gron} implies
\be
M(t)\le C_3 M(0)
\label{2MM0}
\ee
for $0\le t\le t_1$, where $C_3(T,C_0)=\exp(C_1C_2)C_1(1+C_2)$.
Using \eqref{2Mdelay} we infer that
 $M(t_1+2M(t_1))\le C_4M(0)$ with $C_4= C_3(1+4C_*^3)$.
Provided we assume 
\[
M(0)\le\delta_1:=\frac12\frac{\eps_1}{C_4},
\]
it follows that $M(t_1+2M(t_1))<\eps_1$, and since $M$ is continuous,
$\Omega$ is open in $[0,\tilde T]$.

Consequently we have $\tilde T\in\Omega$. Putting $t_1=\tilde T$,
this means we have $M(T)\le\eps_1$ and $M(T)\le C_4M(0)$
if $M(0)\le\delta_1$.
This finishes the proof of Proposition~\ref{P.lip}.

\bigskip\noindent{\bf Proof of Theorem~\ref{T.eu}:}
Uniqueness follows immediately from Proposition~\ref{P.lip}.
To prove existence for arbitrary $v_0\in\RCD$, by Proposition~\ref{P.lip}
and Lemma~\ref{L.theta1} it evidently suffices to prove global existence
for $v_0$ in a dense set of $\RCD$. Solutions in general are constructed
by passing to the limit in $C([0,T],\RCD)$ for every $T>0$.

\begin{lemma} \label{L.dense}
The set of functions in $\RCD$ that take a finite number of values
is dense in $\RCD$.
\end{lemma}

\begin{proof} Let $v_0\in\RCD$ and let $\eps>0$. Let $y_j=\frac12\eps j$
for $j=0,1,\ldots$, and let
\[
v_\eps(\F) = \min\{ y_j\mid y_j\ge v_0(\F)\}
\]
for $\F\in[0,1]$. It is easy to see that $v_\eps$ has a finite number of
values, that $v_\eps\in\RCD$, and
$\|v_\eps-v_0\|<\eps$. This proves the lemma.
\end{proof}

Suppose, then, that $v_0\in\RCD$ takes a finite number of values
$y_0>\ldots>y_N=0$.  Then with $\F_j=\inf\{\F\mid v_0(\F)=y_j\}$,
we have  $0=\F_0<\ldots<\F_N\le1$ and
$v_0(\F)=y_j$ for $\F\in[\F_j,\F_{j+1})$, $j=0,\ldots,N-1$.
We start to construct a solution by solving the system of ordinary
differential equations
\be
w_j'(t) = w_j(t)^{1/3}\Theta(t)-1, \qquad j=0,\ldots,N-1,
\label{2.wj}
\ee
with
\be
\Theta(t) = \F_N\left/ \sum_{j=0}^{N-1} w_j(t)^{1/3}(\F_{j+1}-\F_j) \right.
\label{2.Theta}
\ee
and $w_j(0)=y_j$, on a maximal interval $[0,t_N)$ in which
$\min w_j(t)>0$. The solution is smooth and $w_j(t)>w_{j+1}(t)$
by backwards uniqueness for the equation $w'=w^{1/3}\Theta-1$.
The quantity
\[
\sum_{j=0}^{N-1} w_j(t)(\F_{j+1}-\F_j)
\]
is conserved in time. Without loss of generality we can assume
this quantity is 1.

We can estimate $\Theta(t)\le w_0(t)^{2/3}$ so
$w_0'\le w_0$ and hence $w_0(t)\le e^ty_0$.
If $t_N<\infty$, then, it follows that the smallest component vanishes,
i.e., $w_{N-1}(t_N^-)=\lim_{t\nearrow t_N}w_{N-1}(t)=0$.

For $t\in[0,t_N)$ we define $v(t,\F)=w_j(t)$ for
$\F\in[\F_j,\F_{j+1})$, $j=0,\ldots,N-1$, and let $\theta=\Theta$.
This yields a solution of equations \eqref{2veqn} and \eqref{2vint}
for $t\in[0,t_N)$.
As $t\to t_N$ from below, the limits $v(t_N^-,\F)$ and $\theta(t_N^-)$
exist. The solution can then be re-initialized at time $t_N$ with one
less component ($N$ replaced by $N-1$). After some finite number of such
steps the solution must exist globally.

Thus, for $v_0\in\RCD$ with a finite number of values, a global solution
exists. Theorem~\ref{T.eu} follows.

\bigskip\noindent{\bf Proof of Theorem~\ref{T.lip}:}
The additional restriction imposed in Proposition \ref{P.lip} can be removed
now by considering convex combinations of initial data. Given $T$, $C_0$,
$v_1$, and $v_2$ as stated, let $C>0$, $\delta>0$ be as given by Proposition
\ref{P.lip}, and let $M_0=\|v_1(0,\cdot)-v_2(0,\cdot)\|$. 
Fix an integer $n>M_0/\delta$, and for $j=0,1,\ldots,n$ let
\[
x_j(\F)= \left(1-\frac{j}{n}\right) v_1(0,\F) + \left(\frac{j}{n}\right) v_2(0,\F)
\]
for $\F\in[0,1]$. Then $x_j\in\RCD$, $x_j(0)\le C_0$ for all $j$ and 
$\|x_{j+1}-x_j\|=M_0/n<\delta$. By the existence theorem there exist 
corresponding solutions $v=\tilde v_j$ to \eqref{2vint}--\eqref{2veqn}
with $\tilde v_j(0,\cdot)=x_j$, and by Proposition \ref{P.lip} we have
\[
\sup_{0\le t\le T} \|\tilde v_{j+1}(t,\cdot)-\tilde v_j(t,\cdot)\|\le
C\|x_{j+1}-x_j\| = C M_0/n.
\]
Since $v_1-v_2=\sum_{j=0}^{n-1}(\tilde v_{j+1}-\tilde v_j)$, 
using the triangle inequality we find that 
\[
\sup_{0\le t\le T} \| v_1(t,\cdot)-v_2(t,\cdot)\| \le CM_0,
\]
as desired.
\section{Measure-valued solutions}

Our aim here is to describe a precise correspondence between the solutions
$v(t,\F)$ of Theorem \ref{T.eu} and measure-valued weak solutions
$\nu_t$ of \eqref{1fevol}, and to show that the metric $\|v_1-v_2\|$ on
$\RCD$ corresponds to the $L^\infty$ Wasserstein metric on the space
$\Psp$ of (Borel) probability measures on $[0,\infty)$ with compact support.
Theorem \ref{T.meas} then follows as a corollary of Theorems \ref{T.eu}
and \ref{T.lip}.

We begin with a technical lemma on generalized inverses of 
increasing functions. 

\begin{lemma} \label{L.incr}
Suppose $b>0$ and $w\colon[0,b]\to\R$ is a left continuous increasing
function with $w(0)=0$. Let $b\gi=w(b)$ and define 
$w\gi\colon[0,b\gi]\to\R$ by 
\[
w\gi(y)= \cases{
\sup\{ x\mid w(x)<y\}, & $0<y\le b\gi$, \cr
0, & $y=0$.}
\]
Then $w\gi$ is left continuous and increasing, and moreover, 
\[
w\git=w.
\]
\end{lemma}

\begin{proof} 
Clearly $w\gi$ is increasing. Given $y\in(0,b\gi]$ and $\eps>0$, put
$\bar x=w\gi(y)$ and $2\delta= y-w(\bar x-\eps)$. Then $\delta>0$
and $w(\bar x-\eps)<y-\delta$, hence 
$\bar x-\eps<w\gi(y-\delta)\le \bar x$. It follows $w\gi$ is left continuous.

To show $w\git=w$, it suffices to show that for $0<x<b$,
\[
w(x-\eps)\le w\git(x)\le w(x+\eps)
\]
for all sufficiently small $\eps>0$. Let $\bar y = w\git(x) = 
\sup\{y\mid w\gi(y)<x\}$. Then for all $\eps_0>0$, 
$w\gi(\bar y+\eps_0)\ge x$, hence for any small $\eps>0$ we have
$w(x-\eps)<\bar y+\eps_0$, therefore $w(x-\eps)\le\bar y$.

For the reverse inequality there are two cases: 
If $\bar x=w\gi(\bar y)<x$ then for small
$\eps>0$ we have
\be
\bar y\le w(\bar x+\eps) \le w(x+\eps).
\label{3.wup}
\ee
Otherwise $\bar x\ge x$, and since $w$ is left continuous, $\bar x =x$.
In this case, \eqref{3.wup} again holds. This finishes the proof.
\end{proof}

If $w$ is continuous and strictly increasing, then
$w\gi$ is the inverse function of $w$.

Given a probability measure $\nu$ with compact support 
$[0,\bar v]\subset[0,\infty)$, we associate the distribution function 
$F_\nu\colon[0,\infty)\to[0,1]$  given by
\be
F_\nu(x)=\cases{ \nu([0,x)) & $x>0$\cr 0, & x=0.}
\label{3.Fdef}
\ee
$F_\nu$ is left continuous and increasing, and $F_\nu$ determines $\nu$
(that is, the values of $F_\nu$ determine the values of $\nu$ on all Borel sets).
We associate a decreasing function $v=\vnu(\nu)$ to $\nu$ via
$v(x)=F_\nu\gi(1-x)$ for $x\in[0,1]$. (Here, $F_\nu\gi$ is the generalized
inverse of the restriction of $F_\nu$ to $[0,\bar v+1]$.)
That is,
\be
v(x)= \cases{\sup\{ y\mid F_\nu(y)<1-x\}, & $0\le x<1$,\cr 0,& $x=1$.}
\label{3.vdef}
\ee
With the notation $Rv(x)=v(1-x)$ we have $\vnu(\nu)=R(F_\nu\gi)$. 
The first part of Lemma \ref{L.incr} implies $\vnu(\nu)\in\RCD$, thus the
map $\vnu\colon\Psp\to\RCD$. (Recall $\Psp$ is the set of probability measures
on $[0,\infty)$ with compact support.)

The inverse map to $\vnu$ is given as follows. If $v\in\RCD$ we let
$F=(Rv)\gi$ on $[0,v(0)]$ and put $F(x)=1$ for $x>v(0)$. Then $F$ 
is increasing and left continuous,
and determines a (Borel) probability measure $\nu$ 
for which $F=F_\nu$.  
For later use we note that for any
continuous $f\colon(0,\infty)\to\R$ with compact support, we have
\be
\int_0^1f(v(x))\,dx = \int_0^1f(F\gi(x))\,dx = \int_0^\infty f(y)\,dF(y) = 
\int_0^\infty f(y)\,d\nu(y).
\label{3.chg}
\ee
This follows from \cite[2.5.18(3)]{Fed}, for example. The identity
function $y\mapsto y$ can be approximated uniformly on compact sets
in $[0,\infty)$ by such functions $f$, hence
\be
\int_0^1v(x)\,dx = \int_0^\infty y\,d\nu(y).
\label{3vol}
\ee

We let $\nuv(v)=\nu$, so $\nuv\colon\RCD\to\Psp$.
Lemma \ref{L.incr} implies that we have

\begin{lemma} \label{L.inv}
$\vnu$ and $\nuv$ are inverse maps: $\vnu(\nuv(v))=v$ for all $v\in\RCD$, and
$\nuv(\vnu(\nu))=\nu$ for all $\nu\in\Psp$.
\end{lemma}

We now recall from \cite{GS} that the $L^p$ Wasserstein metric can be defined 
on $\Psp$ as follows. Given $\nu_1$ and $\nu_2$ in $\Psp$,
let $D(\nu_1,\nu_2)$ be the set of probability measures $\mu$ on 
$[0,\infty)\times[0,\infty)$ with marginal distributions $\nu_1$ and $\nu_2$,
that is, for all continuous $\zeta\colon[0,\infty)\to\R$,
\[
\int_0^\infty\int_0^\infty \zeta(x)\,d\mu(x,y)=\int_0^\infty \zeta(x)\,d\nu_1(x)
\]
and
\[
\int_0^\infty\int_0^\infty \zeta(y)\,d\mu(x,y)=\int_0^\infty \zeta(y)\,d\nu_2(y).
\]
If $1\le p<\infty$ then the $L^p$ Wasserstein metric is defined by
\[
d_p(\nu_1,\nu_2) = 
\left(\inf_{\mu\in D(\nu_1,\nu_2)} \int |x-y|^p\,d\mu(x,y) \right)^{1/p}
\]
The $L^\infty$ Wasserstein metric is defined by
\[
d_\infty(\nu_1,\nu_2) = \inf_{\mu\in D(\nu_1,\nu_2)} 
\mbox{$\mu$-$\esssup$} |x-y| .
\]
The measures $\mu$ represent ways to `rearrange  mass' from one distribution
into the other, and the $L^p$ Wasserstein metrics measure the least costly way
to do this according to the notion of cost indicated.

\begin{lemma}
Given $\nu_1$, $\nu_2\in\Psp$, let $v_1=\vnu(\nu_1)$, $v_2=\vnu(\nu_2)$.
Then for $1\le p<\infty$ we have
\[
d_p(\nu_1,\nu_2)= \left(\int_0^1 |v_1(\F)-v_2(\F)|^p\,d\F \right)^{1/p},
\]
and 
\[
d_\infty(\nu_1,\nu_2)=\|v_1-v_2\|.
\]
\end{lemma}

\begin{proof}
The assertion for $1\le p<\infty$ follows from \cite{Rach}, see pp 28--30 and
Corollary 7.3.6, which yields that 
\[
d_p(\nu_1,\nu_2)= \left(\int_0^1 |F_{\nu_1}\gi(\F)-F_{\nu_2}\gi(\F)|^p
\,d\F \right)^{1/p}.
\]
Then Proposition 3 of \cite{GS} asserts that 
$\lim_{p\to\infty} d_p(\nu_1,\nu_2)=d_\infty(\nu_1,\nu_2).$
Since $v_1$ and $v_2$ are right continuous, it follows
\beqna
d_\infty(\nu_1,\nu_2) &=& \lim_{p\to\infty} 
\left(\int_0^1 |v_1(\F)-v_2(\F)|^p\,d\F \right)^{1/p} 
\\&=& \esssup_{[0,1]}|v_1(\F)-v_2(\F)| = \sup_{[0,1]} |v_1(\F)-v_2(\F)|.
\eqna
\end{proof}

\begin{corollary} \label{C.isom}
Let $\Psp$ have the topology induced by $d_\infty$. Then $\Psp$
is complete, and the map $\vnu\colon\Psp\to\RCD$ is an isometric 
isomorphism of complete metric spaces.
\end{corollary}

The completeness of $\Psp$ with respect to the metric $d_\infty$
was established in \cite{GS}.

The correspondence between volume orderings $v\in\RCD$ and volume
distributions $\nu\in\Psp$ has been established. Now we seek to show that
this correspondence maps solutions to weak solutions and vice-versa.

\begin{prop} \label{P.vtonu}
Let $\theta\in L^\infty_\loc(0,\infty)$, $v\in C([0,\infty),\RCD)$ be a
solution of \eqref{2veqn}. 
For each $t\ge0$, let $\nu_t=\nuv(v(t,\cdot))$.
Then $\nu\colon[0,\infty)\to\Psp$ is locally Lipschitz, and
$\nu$ is a weak solution of \eqref{1fevol}.
\end{prop}

\begin{proof}
From \eqref{2veqn} we have that $v\colon[0,\infty)\to\RCD$ is locally
Lipschitz, therefore $\nu\colon[0,\infty)\to\Psp$ is locally Lipschitz
by Corollary \ref{C.isom}. 

Let $\zeta:(0,\infty)\times(0,\infty)\to\R$ be smooth with compact
support. Then for all $\F\in[0,1]$, $t\mapsto \zeta(t,v(t,\F))$ is
Lipschitz continuous and we have 
\beqna
0&=& 
\int_0^\infty \frac{d}{dt}\zeta(t,v(t,\F)) \,dt  
\\&=&  
\int_0^\infty \D_t\zeta(t,v(t,\F)) + 
\vdot(v(t,\F),\theta(t))\D_v\zeta(t,v(t,\F)) \,dt.
\eqna
Let $F(t,\cdot)=(Rv(t,\cdot))\gi$ and let $\F(t,\cdot)=1-F(t,\cdot)$.
We integrate over $\F\in[0,1]$, use Fubini's theorem, and change
variables using \eqref{3.chg}. We obtain
\beqna
0 &=& 
-\int_0^\infty\int_0^\infty\left(\D_t\zeta(t,v) + 
\vdot(v,\theta(t))\D_v\zeta(t,v) \right)\, d\F(t,v)\,dt
\\&=& 
\int_0^\infty\int_0^\infty\left(\D_t\zeta(t,v) + 
\vdot(v,\theta(t))\D_v\zeta(t,v) \right)\, d\nu_t(v)\,dt.
\eqna
Thus $(\theta,\nu)$ is a weak solution in the sense of Theorem \ref{T.meas},
as claimed.
\end{proof}

\begin{prop} \label{P.nutov}
Suppose that $\theta\in L^\infty_\loc(0,\infty)$ and $\nu\colon[0,\infty)\to\Psp$
is locally Lipschitz, and $\nu$ is a weak solution of \eqref{1fevol}.
Let $v(t,\cdot)=\vnu(\nu_t)$ for each $t\ge0$. Then $v$ is a solution
of \eqref{2veqn}.
\end{prop}

\begin{proof}
Given $\theta$ and $v$ as described, the map $v\colon[0,\infty)\to\RCD$ is 
locally Lipschitz. We consider test functions $\zeta$ of the form
\be
\zeta(t,v)=\xi(t)\eta(v)
\label{3.zform}
\ee
where the functions $\xi$, $\eta\colon\R\to\R$ are smooth with
compact support in $(0,\infty)$. Using this form together
with the fact that $(\theta,\nu)$ form a weak solution to \eqref{1fevol},
and using the change of variables from \eqref{3.chg} as previously, we find
that
\[
0=\int_0^\infty\int_0^1
\xi'(t)\eta(v(t,\F)) + \xi(t)\eta'(v(t,\F))\vdot(v(t,\F),\theta(t)) \, d\F\,dt.
\]
Using Fubini's theorem and integrating by parts in time, this gives
\be
0=\int_0^1\int_0^\infty
\xi(t)\tilde\eta(v(t,\F))(\vdot(v(t,\F),\theta(t))-\D_tv(t,\F))\,dt\,d\F.
\label{3.wk1}
\ee
where $\tilde\eta=\eta'$.
This formula is justified since for each $\F\in[0,1]$, $v(\cdot,\F)$ is Lipschitz,
hence differentiable almost everywhere. We note that since $v$ is bounded
on compact sets, $\tilde \eta$ can be chosen to agree on the range of $v$ 
with an arbitrary smooth function with compact support in
$(0,\infty)$. 
We do this and drop the tilde. 
Furthermore, we note that by Lebesgue's dominated convergence theorem,
\eqref{3.wk1} remains valid for any $\xi$ with compact support in 
$(0,\infty)$ that is the bounded pointwise limit
$\xi=\lim_{n\to\infty}\xi_n$ of a sequence of smooth $\xi_n$ with compact support
in $(0,\infty)$, and similarly for $\tilde\eta$. For the moment it will suffice
to consider $\xi$, $\tilde\eta\in \Ccr$, the set of continuous functions
on $(0,\infty)$ with compact support.

For what follows,  we take some care regarding joint measurability in $(t,\F)$
and sets of measure zero.
We fix a representative $\tilde\theta$ in the equivalence class $\theta$,
then drop the tilde. $\D_t v(t,\F)$ need not exist at every point, but equation
\eqref{3.wk1} also holds if $\D_t v$ is replaced by the {upper} derivative
$\vtup$ or the {lower} derivative $\vtdn$, defined by
\[
\vtup(t,\F) = \lim_{\eps\to0}\sup_{0<|h|<\eps} \delvh(t,\F), \quad
\vtdn(t,\F) = \lim_{\eps\to0}\inf_{0<|h|<\eps} \delvh(t,\F),
\]
where
\[
\delvh(t,\F) = \frac{v(t+h,\F)-v(t,\F)}{h}.
\]
Since $v$ is locally Lipschitz in $t$ uniformly in $\F$,
$\vtup$ and $\vtdn$ are bounded on compact sets, and $\vtdn\le\vtup$.

\begin{lemma} \label{L.updn}
As maps from $(0,\infty)\times[0,1]\to\R$, $\vtup$ and $\vtdn$ are Borel 
measurable. Moreover, $\vtup=\vtdn$ almost everywhere in $(0,\infty)\times[0,1]$.
\end{lemma}

\begin{proof}
Since $v$ is continuous in $t$ uniformly in $\F$ and is right continuous and
decreasing in $\F$, $v$ is lower semicontinuous, hence Borel. Suppose 
$0<t_1<t_2<\infty$, then for $0<|h|<t_1$ the map $\delvh$ is Borel on 
$[t_1,t_2]\times[0,1]$. Let $\{h_j\}$ be a dense sequence in $(-1,0)\cup(0,1)$.
Since the maximum of two Borel functions is Borel and pointwise limits of
sequences of Borel functions are Borel, and pointwise we have
\[
\sup_{0<|h|<\eps} \delvh(t,\F) = \sup_{|h_j|<\eps}\delvhj(t,\F)
=\lim_{k\to\infty} \max_{j\le k\atop |h_j|<\eps} \delvhj(t,\F),
\]
by taking $\eps$ to zero along a sequence it follows that $\vtup$ is Borel
on $[t_1,t_2]\times[0,1]$, hence on $(0,\infty)\times[0,1]$.
A similar argument applies for $\vtdn$.

Now we have that the set $Z=\{(t,\F)\mid (\vtup-\vtdn)(t,\F)>0\}$ is a Borel set.
We know that for each $\F$, $v(\cdot,\F)$ is differentiable almost 
everywhere, so $(\vtup-\vtdn)(t,\F)=0$ for almost every $t>0$. 
Fubini's theorem now implies that $Z$ has Lebesgue measure zero in
$(0,\infty)\times[0,1]$.
\end{proof}

Returning to \eqref{3.wk1}, we can now apply Fubini's theorem and deduce that
for almost every $t$, $(\vtup-\vtdn)(t,\F)=0$ for almost every $\F$, and 
with 
\[
J_\eta(t)= \int_0^1\eta(v(t,\F))(\vdot(v(t,\F),\theta(t))-\vtup(t,\F))\,d\F,
\]
we have that for any $\eta\in\Ccr$,
$\int_0^\infty\xi(t)J_\eta(t)\,dt=0$ for all $\xi\in\Ccr$.  Therefore, given $\eta$
there is a set $\Omega_\eta\subset(0,\infty)$ of full measure (meaning the
complement has measure zero), such that $J_\eta(t)=0$ for all $t\in\Omega_\eta$.  

The set $\Ccr$ is separable, so if we do this for a dense sequence $\{\eta_n\}$
we find there is a set $\Omega\subset\cap\Omega_{\eta_n}$ 
of full measure in $(0,\infty)$ such that
for $t\in\Omega$, $J_{\eta_n}(t)=0$ for all $n$. Since any $\eta\in\Ccr$ can be
approximated uniformly by a subsequence of $\{\eta_n\}$, we infer that:

\begin{lemma} \label{L.Jeta}
There is a set $\Omega\subset(0,\infty)$ of full measure,
such that for all $t\in\Omega$, $\D_tv(t,\F)$ exists for almost every $\F\in[0,1]$,
and $J_\eta(t)=0$ for all $\eta\in\Ccr$.
\end{lemma}

\begin{lemma} Let $t\in\Omega$, and suppose $v(t,x)=v(t,y)$ where
$0\le x<y\le1$. Then $\D_tv(t,\F)$ exists for all $\F\in(x,y)$, and 
is constant on this interval.
\end{lemma}

The proof of this lemma is straightforward, using the facts that
$v(t,\cdot)$ is decreasing for every $t$, and $\D_tv(t,\F)$ exists for almost every $\F$.

\begin{lemma} \label{L.zero!}
Let $t\in\Omega$, and let $\bar\F(t)=\sup\{\F\mid v(t,\F)>0\}$. 
Then 
\[
\vdot(v(t,\F),\theta(t))-\D_tv(t,\F)=0
\]
for almost all $\F\in[0,\bar\F(t))$.
\end{lemma}

\begin{proof} We thank B. Kirchheim for the main idea of the following proof. 
Since $t$ is fixed, we suppress indicating dependence on $t$, and
we let $g(\F)=\vdot(v(t,\F),\theta(t))-\vtup(t,\F)$.
We know $g$ is measurable and bounded.
$\F\mapsto v(\F)$ is decreasing, so if $0\le y$ is in the range
of $v$, either the pre-image $v^{-1}(y)$ is a singleton or an interval
of nonzero width. There can be only a countable set of $y$ of the latter
type. Let $\Delta$ be the (countable) set of endpoints of such intervals.
For $x$, $y\in[0,1]\setminus\Delta$, we know that $v(x)=v(y)$ implies
$\vtup(x)=\vtup(y)$, and so $g(x)=g(y)$.

Given any $\eps>0$, let $A_\eps=[0,\bar\F(t))\cap
\{x\in[0,1]\mid g(x)>\eps\}\setminus\Delta$. 
Then $A_\eps$ is measurable, and we claim the measure of $A_\eps$ is 
zero for any $\eps>0$. Suppose not, so $|A_\eps|=2\delta>0$ 
for some $\eps>0$. By Lusin's theorem, 
there is a compact $K\subset A_\eps$ such that $|K|\ge\delta$ and
$v|_K$ is continuous. Then $v(K)$ is compact and is contained in $(0,\infty)$
since $v$ is positive at each point of $[0,\bar\F(t))$.

Apply Lemma \ref{L.Jeta} with $\eta(\F)=\eta_n(\F)=
\max\{0,1-n\dist(\F,v(K))\}$ for $n=1,2,\ldots$ Then $\eta_n$ 
has compact support in $(0,\infty)$ for $n$ sufficiently large and 
converges boundedly pointwise to the characteristic function $\chi_{v(K)}$. 
It follows that 
\[
0=\int_0^1 \chi_{v(K)}(v(\F)) g(\F)\,d\F.
\]
Now if $v(x)\in v(K)$, then $v(x)=v(y)$ for some
$y\in K$, and either $g(x)=g(y)>\eps$ or $x\in\Delta$. It follows
\[
\int_0^1 \chi_{v(K)}(v(\F)) g(\F)\,d\F \ge\eps|K|>0,
\]
yielding a contradiction. Hence $|A_\eps|=0$ for any $\eps>0$.
A similar argument applies for $\{x\mid g(x)<-\eps\}$, and we 
then deduce that $g(\F)=0$ for almost every $\F\in[0,\bar\F(t))$.
This proves the Lemma.
\end{proof}

Note that for $t\in\Omega$ and $\F\in(\bar\F(t),1]$ we have that $v(t,\F)=0$,
and $\D_tv(t,\F)=0$. 

Since $(t,\F)\mapsto v(t,\F)$ is right continuous in $\F$ and locally Lipschitz in 
$t$ uniformly in $\F$, the set 
$Q=\{(t,\F)\mid v(t,\F)>0\}$
is open in $(0,\infty)\times[0,1]$.  Define
\[
g(t,\F)=\vdot(v(t,\F),\theta(t))-\vtup(t,\F),
\]
then $g$ is measurable on $(0,\infty)\times[0,1]$ 
and by Lemma \ref{L.zero!} we have 
\[
\int_0^\infty\int_0^1 \chi_Q|g|\,d\F\,dt = 0.
\]
By Fubini's theorem, we have $\chi_Q g=0$ almost everywhere.
Hence there exists a set $S$ of full measure in $[0,1]$,
such that if $\F\in S$ then $(\chi_Q g)(t,\F)=0$ for almost every $t$,
and $t\mapsto v(t,\F)$ is differentiable almost everywhere.

\begin{lemma}\label{L.eqn}
If $\F\in S$ and $v(t,\F)>0$, then $v(s,\F)>0$
for all $s\in[0,t]$ and 
\be
v(t,\F)=v(0,\F)+\int_0^t \vdot(v(s,\F),\theta(s)) \,ds.
\ee
\end{lemma}

\begin{proof}
For any $t_1\in(0,t)$, if $v(s,\F)>0$ for all $s\in[t_1,t]$
then since $s\mapsto v(s,\F)$ is differentiable 
and $g(s,\F)=0$ almost everywhere in $[t_1,t]$, we have
\[
v(t,\F)=v(t_1,\F)+\int_{t_1}^t \vdot(v(s,\F),\theta(s)) \,ds.
\]
We claim that the set $\{t_1\in[0,t)\mid 
\mbox{ $v(s,\F)>0$ for all $s\in[t_1,t]$}\}$ has the infimum 
$t_*=0$. Note that the set is nonempty 
by the continuity of $s\mapsto v(s,\F)$.
Suppose the infimum $t_*$ is positive. Then 
$v(t_*,\F)=0<v(s,\F)$ for $s\in (t_*,t]$. We know that
$\theta(s)$ is bounded for $s\in[0,t]$, so for some sufficiently
small $h>0$ it follows that $\vdot(v(s,\F),\theta(s))<-1/2$ for
$t_*<s<t_*+h$. Then we have
\[
0<v(t_*+h,\F)=0+\int_{t_*}^{t_*+h} \vdot(v(s,\F),\theta(s))\,ds
< -h/2,
\]
a contradiction. Hence our claim holds: $v(s,F)>0$ for
all $s\in[0,t]$, and the formula asserted in the Lemma follows.
\end{proof}

Now we can finish the proof of Proposition \ref{P.nutov}.
Suppose $t>0$, $\F\in[0,1]$ are arbitrary and $v(t,\F)>0$.
Since $v(t,\cdot)$ is right continuous and decreasing, there exists
a sequence of numbers $\F_n\in S$ such that $\F_n>\F$, 
$\F_n\to\F$ as $n\to\infty$ and $v(t,\F_n)>0$. Using Lemma \ref{L.eqn}
it follows that $v(s,\F)>0$ for all $s\in[0,t]$ and
\beqna
&&v(t,\F)-v(0,\F)-\int_0^t \vdot(v(s,\F),\theta(s))\,ds
\\&& \qquad=\ 
\lim_{n\to\infty} \left(
v(t,\F_n)-v(0,\F_n)-\int_0^t \vdot(v(s,\F_n),\theta(s))\,ds
\right) =0.
\eqna
This completes the proof of 
Proposition \ref{P.nutov}.
\end{proof}

The results asserted in Theorem \ref{T.meas} now follow directly from
Theorems \ref{T.eu} and \ref{T.lip}, with the help of 
Propositions \ref{P.vtonu} and \ref{P.nutov}, Corollary \ref{C.isom},
and equation \eqref{3vol}.
\section{A different conserved quantity}

In the theory of Ostwald ripening, one also encounters an
alternative to the condition that the total particle volume is
conserved in time. If mass in the diffusion field is taken into
account, one finds that a quantity of the form
\be
Q= a\theta(t)+\int_0^1 v(t,\F)\,d\F
\label{4con}
\ee
is conserved instead, where $a>0$ is a constant. 
$\theta(t)$ need no longer be positive.

In terms of the theory developed in this paper, the constraint
\eqref{4con} is simpler to deal with that the constraint of
constant volume. One has the bound
\[
\theta(t) \le Q/a,
\]
and when comparing two solutions of \eqref{2veqn}, 
one can use the arguments of Lemma \ref{L.Mest} and replace
the use of Lemma \ref{L.thetaest} by the simpler estimate
\be
|\theta_1(t)-\theta_2(t)| \le a^{-1}\|v_1(t,\cdot)-v_2(t,\cdot)\| +
a^{-1}|Q_1-Q_2|.
\ee
From the standard Gronwall's inequality, one easily deduces the
a priori estimate asserted in the following result. The existence
and uniqueness proofs go the same as in section 3.

\begin{theorem} \label{T.ivp2}
Let $v_0\in\RCD$, $Q\in\R$. Then there exists a unique function
$v\in C([0,\infty),\RCD)$ such that, with $\theta(t)$ determined by
\eqref{4con}, we have
\[
v(t,\F)=v_0(\F)+\int_0^t(v(s,\F)^{1/3}\theta(s)-1)\,ds
\]
whenever $v(t,\F)>0$.

Given $T>0$, $C_0>0$, there exists a positive constant $C$ 
such that, given two solutions as above which also satisfy
$\max(Q_1,Q_2)\le C_0$,
then
\[
\sup_{0\le t\le T} \|v_1(t,\cdot)-v_2(t,\cdot)\| \le 
C \left(\|v_1(0,\cdot)-v_2(0,\cdot)\|  +|Q_1-Q_2| \right).
\]
\end{theorem}

Using the correspondence $v(t,\cdot)\mapsto \nu_t=\nuv(v(t,\cdot))$
and its inverse as in section 3, from Propositions \ref{P.vtonu} and \ref{P.nutov}
one may deduce directly the following corollary of Theorem \ref{T.ivp2}.

\begin{theorem} \label{T.meas2}
Given $\nu_0\in\Psp$, $Q\in\R$, there exists a unique map
$t\mapsto\nu_t$ that is locally Lipschitz from $[0,\infty)$ into $\Psp$
such that, with $\theta(t)$ determined by the relation
\[
Q= a\theta(t)+\int_0^\infty v\,d\nu_t(v),
\]
we have
\[
\int_0^\infty\int_0^\infty 
\D_t\zeta(t,v)+\vdot(v,\theta(t))\D_v\zeta(t,v)\ d\nu_t(v)\,dt = 0
\]
for all smooth $\zeta\colon(0,\infty)\times(0,\infty)\to\R$ with compact support.

Given any $T>0$, $C_0>0$, there exists $C>0$ such that,
if two such weak solutions $\nuone$, $\nutwo$ are
given, which satisfy $\max(Q_1,Q_2)\le C_0$,
then
\[
\sup_{0\le t\le T} d_\infty(\nuone_t,\nutwo_t) \le 
C\left( d_\infty(\nuone_0,\nutwo_0) +|Q_1-Q_2|\right).
\]

\end{theorem}

\bigskip\noindent{\bf Acknowledgments.}
The authors gratefully acknowledge discussions
with S. M\"uller, B. Kirchheim, F. Otto, S. Luckhaus, 
G. Friesecke, and R. Kohn. 
Part of this work was performed at the
Max-Planck Institute for Mathematics in the Sciences in Leipzig.
This work was partly supported by the National Science Foundation
under grant DMS-9704924, and the 
Deutsche Forschungsgemeinschaft through the SFB 256.


\begin{thebibliography}{99999}
\frenchspacing
\bibitem{Fed} H. Federer,
{\it Geometric Measure Theory}, Springer-Verlag, New York, 1969.

\bibitem{GS} C. R. Givens and R. M. Shortt,
A class of Wasserstein metrics for probability distributions,
{Mich.\ Math.\ J. \bf 31} (1984) 231--240.

\bibitem{LS} I. M. Lifshitz and V. V. Slyozov,
The kinetics of precipitation from supersaturated solid solutions,
{J.\ Phys.\ Chem.\ Solids \bf 19} (1961) 35--50.


\bibitem{Rach} S. T. Rachev,
{\it Probability Metrics and the Stability of Stochastic Models},
Wiley, New York, 1991.

\bibitem{Wa} C. Wagner,
Theorie der Alterung von Niederschl\"agen durch Uml\"osen,
{Z. Elektrochem. \bf 65} (1961) 581--594.

\end {thebibliography}
\end{document}